\documentclass[12pt,twoside]{article}
\usepackage{amsmath,amssymb}
\usepackage{graphicx}
\usepackage{epstopdf}
\usepackage{float}
\usepackage{mathrsfs,amssymb}

\allowdisplaybreaks[4]

\title{\Large
Single-valley-extended solutions with platforms of FKS equation
\thanks{A Project Supported by Scientific Research Fund of SiChuan Provincial Education Department (18ZA0274).}
}

\author{\small
{\sc  Yong-Guo SHI}\footnote{Email: scumat@163.com}
\\
{\small
Numerical Simulation Key Laboratory of Sichuan Province,}
\\
{\small  College of Mathematics and Information Science,
Neijiang Normal University}\\
{\small Neijiang, Sichuan 641100, P. R. China  }
  }
\date{}

\begin{document}
\maketitle
\renewcommand{\theequation}{\thesection.\arabic{equation}}
\newtheorem{lm}{Lemma}[section]
\newtheorem{defi}{Definition}[section]
\newtheorem{prob}{Problem}[section]
\newtheorem{thm}{Theorem}[section]
\newtheorem{pro}{Proposition}[section]
\newtheorem{exmp}{Example}[section]
\newtheorem{rmk}{Remark}[section]
\newtheorem{alg}{Algorithm}[section]
\newtheorem{cor}{Corollary}[section]
\newtheorem{con}{Conjecture}[section]
\renewcommand\figurename{\rm Fig.}
\renewcommand\tablename{\rm Tab.}

\begin{abstract}
The Feigenbaum-Kadanoff-Shenker (FKS) equation is a nonlinear iterative functional equation, which
characterizes the quasiperiodic route to chaos for circle maps. Instead of
FKS equation, we put forward the second type of FKS equation.
We use the iterative construction method to obtain a new kind of continuous solutions with platforms for this functional equation.
We show some properties of these continuous solutions, and give one relationship of some solutions between FKS equation and the second type of FKS equation.
\end{abstract}

\baselineskip 15pt   
\parskip 10pt         

\thispagestyle{empty}
\setcounter{page}{1}


\section{Introduction}

In 1982, Feigenbaum, Kadanoff, and Shenker (\cite{FKS82PD}) analyzed the quasiperiodic route to chaos for maps of the circle,
and formulated the following equation
\begin{equation}\label{FKS}
\left\{ {\begin{array}{*{20}l}
   {g(g(\epsilon^2 x))=-\epsilon g(x),} &
\\
    {g(0)=1,} &
\end{array}} \right.
\end{equation}
where $\epsilon\in (0,1)$  and $g:[-1,1]\to [-1,1]$ is an unknown even unimodal mapping.
We refer to this equation as the {\it Feigenbaum-Kadanoff-Shenker (FKS) equation} (\cite{BrigDixonSzekeres1998IJBC}).
This is a kind of nonlinear iterative functional equation. In general, it is difficult to solve this kind of equations.

For Eq.(\ref{FKS}), as the hypothesis in \cite{FKS82PD,BrigDixonSzekeres1998IJBC,Brig2011}, $g$ decreases
monotonically from $g(0)=1$ to its zero $g(\alpha)=0$ and then to
$g(1)=-\epsilon$.
Mestel in \cite{Mestel} considered analytic solutions of Eq.(\ref{FKS}),
(also cf. \cite{EckmannEpstein1986}).
Recently, we give continuously differentiable unimodal solutions in \cite{Shi2020DCDSB}.
We put forward the second type of FKS equation
\begin{equation}\label{FKS2}
\left\{ {\begin{array}{*{20}l}
   {f(f(\epsilon^2 x))=\epsilon f(x),} &
\\
    {f(0)=1,} &
\end{array}}
\right.
\end{equation}
where $\epsilon\in (0,1)$  and $f:[0,1]\to [0,1]$ is an unknown function.
In \cite{Shi2019AM}, we show that finding single-peak even solutions of FKS equation (\ref{FKS})
is equivalent to finding single-valley solutions of the second type of FKS equation (\ref{FKS2}).
We construct all single-valley solutions of (\ref{FKS2}).

A natural question that arises here is whether the equation (\ref{FKS2}) has other kinds of continuous solutions
as the Feigenbaum equation (cf. \cite{XS1997JM})?

This paper is orgnized as follows. Section 2 are devoted to some properties of continuous solutions of (\ref{FKS2}).
In section 3, we give the definition of single-valley-extended solutions with platforms  of (\ref{FKS2}), and show some properties of this new kind of solutions. In section 4, we construct all single-valley-extended solutions with platforms for (\ref{FKS2}).
In the last section, we give one relationship between single-valley-extended solutions with platforms of Eq.(\ref{FKS2}) and even solutions of Eq.(\ref{FKS}).


\section{Properties of continuous solutions}


In this section, we will present some properties of solutions of (\ref{FKS2}).
Let $f$ be a continuous solution of (\ref{FKS2}). It is clear that $f(1) = \epsilon$, $f^2(\epsilon^2)=\epsilon^2$.
Furthermore, we have the following Lemmas.
\begin{lm}\label{P1} Let $\alpha$ be a minimum point of $f$. Then $f(\alpha)=0$, and $\alpha>\epsilon^2$.
\end{lm}
\par
\emph{Proof}.
By (\ref{FKS2}), we have
$$
f(\alpha)=\epsilon^{-1}f^2(\epsilon^2\alpha)\geq \epsilon^{-1}f(\alpha).
$$
Thus $f(\alpha)=0$. Assume that $\alpha \leq \epsilon^2$. Then $\frac{\alpha}{\epsilon^2}\in [0,1]$.
Since
$$
\epsilon f\left(\frac{\alpha}{\epsilon^{2}}\right)=f^{2}(\alpha)=f(0)=1,
$$
one has $f(\frac{\alpha}{\epsilon^2})=\epsilon^{-1}>1$. This is a contradiction.
Therefore $\alpha>\epsilon^2$.
\hfill$\Box$

From  Lemma \ref{P1}, we see that all minimum points of $f$ lie in $(\epsilon^2,1)$.

\begin{lm}\label{P2}
Let $\alpha$ be a minimum point of $f$. If
there is a point $x\in [0,\epsilon^2]$ such that $f(x)=\alpha$, then there exists a minimum point $\beta$ of $f$ such that
$x=\epsilon^2 \beta$.
\end{lm}
\par
\emph{Proof}.
It follows from $f(x)=\alpha$ that $f(f(x))=f(\alpha)=0$. By  (\ref{FKS2}) and $x\in [0,\epsilon^2]$, we have
$$
f\left(\frac{x}{\epsilon^{2}}\right)=\epsilon^{-1} f^{2}(x)=0,
$$
which implies that $\frac{x}{\epsilon^2}$ is a minimum point of $f$, denoted by $\beta$. Therefore $x=\epsilon^2 \beta$.
\hfill$\Box$

\begin{lm}\label{P3}
Let $S$ be the set that contains all the minimum points of $f$. Then $S$ is a singleton or a closed subinterval in $(\epsilon^2,1)$.
\end{lm}
\par
\emph{Proof}. By \cite{Shi2019AM}, we see that $S$ is a singleton. Assume that $S$ is not a singleton in $(\epsilon^2,1)$. One can see
$S\neq \emptyset$. Let $u=\inf S$ and $v=\sup S$. By the continuity of $F$, $u\neq v$. By Lemma \ref{P2}, there exist $\beta_1,\beta_2\in S$
such that $f(\epsilon^2 \beta_1)=u$ and $f(\epsilon^2 \beta_2)=v$. Let $\eta$ be a maximum point of $f$ in $[u,v]$. Then
there exists a point $\xi$ between $\epsilon^2\beta_1$ and $\epsilon^2\beta_2$  such that $f(\xi)=\eta$.
Thus $\xi/\epsilon^2 \in(\beta_1,\beta_2) \subset [u,v]$. Consequently, we have
$$
f\left(\frac{\xi}{\epsilon^{2}}\right)=\epsilon^{-1} f^{2}(\xi)=\epsilon^{-1} f(\eta) \leq f(\eta).
$$
Since $0<\epsilon<1$, we have $f(\eta)=0$.  Therefore every point in $[u,v]$ is a minimum point of $f$, and $S=[u,v]$.
\hfill$\Box$

\begin{lm}\label{P4}
Let $[u,v]$ be the interval consisting of all minimum points of $f$. Then the following facts hold:
\begin{itemize}
  \item[{\rm (i)}]  $f([0,\epsilon^2 u])=[v,1]$,  $f(\epsilon^2 u)=v$ and $f([\epsilon^2 u,\epsilon^2 v])\subseteq [u,v]$.
    \item[{\rm (ii)}] $f([v,1])=[0,\epsilon]$.
      \item[{\rm (iii)}] $f(\epsilon^2)>v$ or $f(\epsilon^2)<u$.
    \item[{\rm (iv)}] If $f(\epsilon^2)>v$, then
    $$
   f([\epsilon^2 v,\epsilon^2])\subseteq [v,1], f(\epsilon^2 v)=v.
    $$
  \item[{\rm (v)}] If $f(\epsilon^2)<u$, then
    $$
  f([\epsilon^2 u,\epsilon^2 v])=[u,v],   f([\epsilon^2 v,\epsilon^2])\subset[\epsilon^2u,u], f(\epsilon^2 v)=u.
    $$
\end{itemize}
\end{lm}
\par
\emph{Proof}. For fact (i), for $\forall x\in [0, \epsilon^2 u)$, one can see that
$f(f(x))= \epsilon f(\frac{x}{\epsilon^2})>0$. Thus $f(x) \notin [u,v]$.
Since $f(0)=1$, by the continuity of $f$ and Lemma \ref{P2}, we have
$$
f([0,\epsilon^2 u])=[v,1], f(\epsilon^2 u)=v.
$$
For $\forall x\in [\epsilon^2 u, \epsilon^2 v]$, we have
$x/\epsilon^2\in [u,v]$. Thus $f(f(x))= \epsilon f(\frac{x}{\epsilon^2})=0$.
Consequently, $f(x) \in [u,v]$, i.e.,
$$
f([\epsilon^2 u,\epsilon^2 v])\subseteq [u,v].
$$

For fact (ii),  let $\eta$ be a maximum point of $f$ on the subinterval $[v,1]$.
By fact (i),
there exists a point $\xi\in [0,\epsilon^2 u]$  such that $f(\xi)=\eta$. Thus
$$
f(\eta)=f(f(\xi))=\epsilon f\left(\frac{\xi}{\epsilon^{2}}\right) \leq \epsilon.
$$
Together with $f(v)=0$ and $f(1)=\epsilon$, we have $f([v,1])=[0,\epsilon]$.

For fact (iii), assume that $v \leq f(\epsilon^2)\leq u$. Thus
$f(f(\epsilon^2))=0$, which is a contradiction with $f^2(\epsilon^2)=\epsilon^2$.

For fact (iv), for $\forall x\in (\epsilon^2 v, \epsilon^2]$, we have
$x/\epsilon^2\in (v,1]$. Thus $f(f(x))= \epsilon f(\frac{x}{\epsilon^2})>0$.
Thus $f(x) \notin [u,v]$. Since $f(\epsilon^2)>v$, by the continuity of $f$ and Lemma \ref{P2}, we have
$$
f([\epsilon^2 v,\epsilon^2])\subseteq [v,1].
$$
By lemma \ref{P2}, we have $f(\epsilon^2 v)=v$.

For fact (v), for $\forall x\in (\epsilon^2 v, \epsilon^2]$, we have
$x/\epsilon^2\in (v,1]$. Thus $f(f(x))= \epsilon f(\frac{x}{\epsilon^2})>0$.
Thus $f(x) \notin [u,v]$. Since $f(\epsilon^2)<u$, by the continuity of $f$, we have
$$
f([\epsilon^2 v,\epsilon^2])\subseteq [0,u].
$$
By lemma \ref{P2}, we have $f(\epsilon^2 v)=u$.

For $\forall x\in [\epsilon^2 v, \epsilon^2]$, we have $x/\epsilon^2\in [v,1]$. It follows from fact (ii) that
$f(x/\epsilon^2) \leq \epsilon$. Together with Lemma \ref{P1}
$$
f(f(x))=\epsilon f\left(\frac{x}{\epsilon^{2}}\right) \leq \epsilon^{2}<v.
$$
Since $f([0,\epsilon^2 u])=[v,1]$ (by fact (i)), we have $f(x)>\epsilon^2 u$.
Thus $f([\epsilon^2 v,\epsilon^2])\subset[\epsilon^2u,u]$ and $f(\epsilon^2)>\epsilon^2 u$.
\hfill$\Box$

\section{Properties of single-valley-extended solutions with platforms}


In this section, we first introduce some definitions.

\begin{defi}\label{SVP}
{\rm A function $f$ on the interval $[a,b]$ is called a {\it single-valley function with a platform}
if there exists a subinterval $[u,v] \subset [a,b]$ such that $f$ is strictly decreasing
on $[a,u]$, $f$ is a constant on $[u,v]$, and $f$ is strictly increasing on $[v,b]$.
}
\end{defi}

\begin{defi}\label{SPS}
{\rm A solution $f$ is called a {\it single-valley-extended
solution  with platforms} of (\ref{FKS2}) if $f$ is a continuous solution of (\ref{FKS2}), and $f$ is a {\it single-valley function with a platform} on $[\epsilon^2, 1]$.
}
\end{defi}

Let $f$ be a single-valley-extended solution with platforms of (\ref{FKS2}).
Let $[u,v]$ be the interval consisting of all minimum points of $f$. Then $\epsilon^2<u<v<1$
and for $\forall x\in [u,v]$, $f(x)=0$. Furthermore, we have the following properties of this kind of solutions.

\begin{lm}\label{P5}
The equation $f(x)=\epsilon x$ has the unique solution $x=1$ on $[u,1]$.
\end{lm}
\par
\emph{Proof}.
It is clear that $x=1$ is a solution of $f(x)=\epsilon x$. Assume that $\eta\in [v,1]$ is another
solution of  $f(x)=\epsilon x$.
From fact (i) of Lemma \ref{P4}, there exists a point $\xi\in [0,\epsilon^2u]$ such that $f(\xi)=\eta$.
Then
$$
f(f(\epsilon^2 \xi))=\epsilon f(\xi))=\epsilon \eta=f(\eta).
$$
It follows from fact (i) of Lemma \ref{P4} that $f(\epsilon^2 \xi)>v$.
Since $f$ is strictly increasing on $[v,1]$, we have
$$
f(\epsilon^2 \xi)=\eta.
$$
By induction, one has $f(\epsilon^{2n} \xi)=\eta$ for $n=1,2,...$.
When $n$ tends to $+\infty$, it follows from the continuity of $f$, we have
$$
\eta=\lim\limits_{n\to +\infty}f(\epsilon^{2n}\xi)=f(0)=1.
$$
Since $f(x)=0$ for $x\in [u,v]$, we have $f(x)=\epsilon x$ has the unique solution $x=1$ on $[u,1]$.
\hfill$\Box$

\begin{lm}\label{P8} (i)
$f(\epsilon^2)\geq \epsilon^2$.
(ii) If $f(\epsilon^2)<u$, then $f([\epsilon^2 v,\epsilon^2])\subseteq[\epsilon^2,u]$.
\end{lm}
\par
\emph{Proof}.
We prove (i) by contradiction. Assume that $s:=f(\epsilon^2)<\epsilon^2$. By (v) of Lemma \ref{P4}, one can see $s>\epsilon^2 u$. It follows from  $f^2(\epsilon^2)=\epsilon^2$ that
$$
f\left(\frac{s}{\epsilon^2} \right)=\epsilon^{-1}f^2(s)=\epsilon^{-1}f(\epsilon^2)=\epsilon\cdot \frac{s}{\epsilon^2}.
$$
It shows that the equation $f(x)=\epsilon x$ has a solution $x=s/\epsilon^2\in (u,1)$. This is a contradiction to Lemma \ref{P5}.

For fact (ii), for $\forall x\in [\epsilon^2 v, \epsilon^2]$, we have $x/\epsilon^2\in [v,1]$. Then
$$
f(f(x))=\epsilon f\left(\frac{x}{\epsilon^{2}}\right).
$$
According to the definition of $f$, the function $F(x)=\epsilon f\left(\frac{x}{\epsilon^{2}}\right)$ is continuous and strictly increasing.
By Lemma 15.1 in \cite[p.297]{Kuczma1968}, we see that $f$ is strictly monotonic on $[\epsilon^2 v, \epsilon^2]$.
If $f(\epsilon^2)<u$, then by (v) of Lemma \ref{P4} we see  $f([\epsilon^2 v,\epsilon^2])\subseteq[\epsilon^2u,u]$ and $f(\epsilon^2 v)=u$.
Consequently, $f$ is strictly decreasing on $[\epsilon^2 v, \epsilon^2]$. The result (ii) follows from $f(\epsilon^2)\geq \epsilon^2$.
\hfill$\Box$

\begin{lm}\label{P6}
If $f(\epsilon^2)>v$, Then for each integer $n\geq 0$,  $f$ is strictly increasing
on $[\epsilon^{2n}v,\epsilon^{2n}]$,  and strictly decreasing on $[\epsilon^{2n+2},\epsilon^{2n}u]$.
\end{lm}
\par
\emph{Proof}.
When $n=0$, the result is obvious.

When $n=1$, let $p$ and $q$ be the restrictions of $f$ to the subintervals $[v,1]$ and $[\epsilon^2,u]$, respectively.
By Lemma \ref{P4}, one can see that
$$
f([\epsilon^4,\epsilon^2 u]\cup [\epsilon^2v,\epsilon^2])\subseteq [v,1].
$$
For $x\in [\epsilon^4,\epsilon^2 u]$, we have $x/\epsilon^2\in [\epsilon^2, u]$. Thus $p(f(x))=\epsilon q(x/\epsilon^2)$. Consequently,
$$
f(x)=p^{-1}(\epsilon q(x/\epsilon^2)).
$$
Since $p^{-1}$ is strictly increasing and $q$ strictly decreasing, we see that
$f$ is strictly decreasing on $[\epsilon^4,\epsilon^2 u]$.

For $x\in [\epsilon^2v,\epsilon^2]$, we have $x/\epsilon^2\in [v, 1]$. Thus $p(f(x))=\epsilon p(x/\epsilon^2)$. Consequently,
$$
f(x)=p^{-1}(\epsilon p(x/\epsilon^2)).
$$
Since $p^{-1}$ and $p$ are strictly increasing, we see that
$f$ is strictly increasing on $[\epsilon^2v,\epsilon^2]$

Then the result follows by induction.
\hfill$\Box$

\begin{lm}\label{P7}
If $f(\epsilon^2)<u$, then for each integer $n\geq 0$,  $f$ is strictly decreasing
on $[\epsilon^{2n}v,\epsilon^{2n}]$ and $[\epsilon^{2n+2},\epsilon^{2n}u]$.
\end{lm}
\par
\emph{Proof}.
When $n=0$, the result is obvious.

When $n=1$, let $p$ and $q$ be the restrictions of $f$ to the subintervals $[v,1]$ and $[\epsilon^2,u]$, respectively.
By Lemma \ref{P8}, one can see that
$$
f([\epsilon^2 v,\epsilon^2])\subset[\epsilon^2,u].
$$
For $x\in [\epsilon^2v ,\epsilon^2]$, we have $x/\epsilon^2\in [v, 1]$. Thus $q(f(x))=\epsilon p(x/\epsilon^2)$. Consequently,
$$
f(x)=q^{-1}(\epsilon p(x/\epsilon^2)).
$$
Since $q^{-1}$ is strictly decreasing and $p$ strictly increasing, we see that
$f$ is strictly decreasing on $[\epsilon^2v ,\epsilon^2]$.

By Lemma \ref{P4}, one can see that
$$
f([\epsilon^4 ,\epsilon^2u])\subset[v,1].
$$
For $x\in [\epsilon^4 ,\epsilon^2u]$,
we have $x/\epsilon^2\in [\epsilon^2, u]$. Thus $p(f(x))=\epsilon q(x/\epsilon^2)$. Consequently,
$$
f(x)=p^{-1}(\epsilon q(x/\epsilon^2)).
$$
Since $p^{-1}$ is strictly increasing, and $q$ strictly decreasing, we see that
$f$ is strictly decreasing on $[\epsilon^4,\epsilon^2 u]$.

Then the result follows by induction.
\hfill$\Box$

\section{Construction of single-valley-extended solutions with platforms}


In this section, we will construct all single-valley-extended solutions with platforms of Eq.(\ref{FKS2}).
According to Lemma \ref{P4}, it suffices to consider the following two cases:
$\epsilon^2 \leq f(\epsilon^2)< u$ and $v< f(\epsilon^2)$.

\subsection{The case $\epsilon^2 \leq f(\epsilon^2)< u$ }

\begin{thm}\label{T1}    Choose arbitrarily $u,v$ and $a$ such that $\epsilon^2 \leq a< u< v<1$.
Define any strictly decreasing continuous function $q: [\epsilon^2,u]\to [0, a]$ and any strictly increasing continuous function
$p: [v,1]\to [0, \epsilon]$ satisfying:
 \begin{itemize}
 \item[(i)] $q(\epsilon^2)=a$,  $q(a)=\epsilon^2$, $p(v)=q(u)=0$, $q(1)=\epsilon$,
   \item[(ii)]  $p(x)=\epsilon x$ has the unique solution $x=1$ on $[v,1]$.
\end{itemize}
Define
\begin{eqnarray}\label{phi0}
\varphi_0(x)=\left\{ {\begin{array}{*{20}l}
    {q(x),} &   x\in [\epsilon^{2},u],\\
        {0,} &   x\in [u,v],\\
        {p(x),} &   x\in [v,1],
\end{array}} \right.
\end{eqnarray}
\begin{eqnarray}\label{vphi1}
\varphi_1(x)=\left\{ {\begin{array}{*{20}l}
    {q^{-1}(\epsilon  p(x/\epsilon^{2})),} &   x\in [\epsilon^{2}v,\epsilon^{2}],\\
    {h(x),} &   x\in [\epsilon^{2}u,\epsilon^{2}v],\\
    {p^{-1}(\epsilon  q(x/\epsilon^{2})),} & x\in [\epsilon^{4},\epsilon^{2}u],
\end{array}} \right.
\end{eqnarray}
where $h(x)$ is any continuous function satisfying $h(\epsilon^{2}v)=u$ and $h(\epsilon^{2}u)=v$.
For $k\geq2$, define
\begin{equation}\label{recurrence2}
\varphi_{k}(x)=p^{-1}(\epsilon  \varphi_{k-1}(x/\epsilon^{2})),\,\, \mbox{for}\,\, x\in [\epsilon^{2k+2},\epsilon^{2k}],
\end{equation}
and
\begin{eqnarray*}
\varphi(x)=\left\{ {\begin{array}{*{20}l}
    {1,} &  x=0,\\
    {\varphi_k(x),} & x\in [\epsilon^{2k+2},\epsilon^{2k}],\,\, k=0,1,2,....
\end{array}} \right.
\end{eqnarray*}
Then $\varphi$ is a single-valley-extended solutions with platforms of Eq.(\ref{FKS2}). Conversely, if
$\varphi_k$ is the restriction on $[\epsilon^{2k+2},\epsilon^{2k}]$ of a single-valley-extended solutions with platforms of Eq.(\ref{FKS2}),
then (i), (ii), (\ref{vphi1}) and (\ref{recurrence2}) hold.
\end{thm}
\par
\emph{Proof}.
One can check that $\varphi_k$ is well-defined, and continuous on $[\epsilon^{2k+2},\epsilon^{2k}]$ for each $k=1,2,...$.
It follows from (\ref{recurrence2}) that
$$
\varphi_{1}(\epsilon^{2}-0)=p^{-1}(\epsilon  \varphi_{0}(1-0))=p^{-1}(\epsilon^2)=\varphi_{0}(\epsilon^{2}).
$$
By induction, we have
$$
\varphi_{k}(\epsilon^{2k}-0)=\varphi_{k-1}(\epsilon^{2k}+0), \,\, k=1,2,....
$$
Thus $\varphi$ is continuous on $(0,1]$.

Now it suffices to show that $\varphi$ is continuous at $x=0$.
Let $\zeta$ be a minimum point of $h$ on $[\epsilon^{2}u,\epsilon^{2}v]$.
Consider the sequence $\{\varphi_{k}(\epsilon^{2k-2}\zeta)\}$.
According to the definition of $\varphi_{k}$, $\{\varphi_{k}(\epsilon^{2k-2}\zeta)\}_{k=2}^{+\infty}$ is a strictly increasing sequence
on $[v,1]$. It follows from the monotone bounded theorem that this sequence has a limit, denoted by  $\xi=\lim\limits_{k\to +\infty}\varphi_{k}(\epsilon^{2k-2}\zeta)$. Then $\xi\in [v,1]$. Thus for $k\geq2$
$$
p(\varphi_{k}(\epsilon^{2k-2}\zeta))=\epsilon\varphi_{k-1}(\epsilon^{2k-4}\zeta).
$$
When $k$ tends to $+\infty$, we have $p(\xi)=\epsilon \xi$. By (ii), one can see that $\xi=1=\varphi(0)$.
Note that $\varphi_{k}$ takes the minimum at $\epsilon^{2k-2}\zeta$. Thus $\varphi|_{[0,\epsilon^{2k-2}\zeta]}$ takes the minimum at $\epsilon^{2k-2}\zeta$, which guarantees the continuity of $\varphi$ at $x=0$.

One can check directly that $\varphi$ is a single-valley-extended solutions with platforms of Eq.(\ref{FKS2}).

Conversely, suppose $\varphi_k$ is the restriction on $[\epsilon^{2k+2},\epsilon^{2k}]$ of a single-valley solution of Eq.(\ref{FKS}).
It follows from Lemma \ref{P1} that $\varphi_0$ satisfies (i), (ii), (\ref{vphi1}) and (\ref{recurrence2}).
\hfill$\Box$

\subsection{The case $f(\epsilon^2)>v$}

\begin{thm} \label{T2}
  Choose arbitrarily $u,v$ and $a$ such that $\epsilon^2 < a< u< v<1$.
Define any strictly decreasing continuous function $q: [\epsilon^2,u]\to [0, a]$ and any strictly increasing continuous function
$p: [v,1]\to [0, \epsilon]$ satisfying:
 \begin{itemize}
 \item[(i)] $q(\epsilon^2)=a$,  $q(a)=\epsilon^2$, $p(v)=q(u)=0$, $q(1)=\epsilon$,
   \item[(ii)]  $p(x)=\epsilon x$ has the unique solution $x=1$ on $[v,1]$.
\end{itemize}
Define
\begin{eqnarray}\label{f0}
\varphi_0(x)=\left\{ {\begin{array}{*{20}l}
    {q(x),} &   x\in [\epsilon^{2},u],\\
        {0,} &   x\in [u,v],\\
        {p(x),} &   x\in [v,1],
\end{array}} \right.
\end{eqnarray}
and
\begin{eqnarray}\label{h1}
\varphi_1(x)=\left\{ {\begin{array}{*{20}l}
    {p^{-1}(\epsilon  p(x/\epsilon^{2})),} &   x\in [\epsilon^{2}v,\epsilon^{2}],\\
    {h(x),} &   x\in [\epsilon^{2}u,\epsilon^{2}v],\\
    {p^{-1}(\epsilon  q(x/\epsilon^{2})),} & x\in [\epsilon^{4},\epsilon^{2}u],
\end{array}} \right.
\end{eqnarray}
where $h(x)$ is any continuous function satisfying $h(\epsilon^{2}v)=h(\epsilon^{2}u)=v$.
For $k\geq2$, define
\begin{equation}\label{h2}
\varphi_{k}(x)=p^{-1}(\epsilon  \varphi_{k-1}(x/\epsilon^{2})),\,\, \mbox{for}\,\, x\in [\epsilon^{2k+2},\epsilon^{2k}],
\end{equation}
and
\begin{eqnarray*}
\varphi(x)=\left\{ {\begin{array}{*{20}l}
    {1,} &  x=0,\\
    {\varphi_k(x),} & x\in [\epsilon^{2k+2},\epsilon^{2k}],\,\, k=0,1,2,....
\end{array}} \right.
\end{eqnarray*}
Then $\varphi$ is a single-valley-extended solutions with platforms  of Eq.(\ref{FKS2}). Conversely, if
$\varphi_k$ is the restriction on $[\epsilon^{2k+2},\epsilon^{2k}]$ of a single-valley-extended solutions with platforms of Eq.(\ref{FKS2}),
then (i), (ii), (iii), (\ref{h1}) and (\ref{h2}) hold.
\end{thm}
\par
The proof is similar to the one of Theorem \ref{T1}. We omit it here.

\section{Relationship}


The following result can provide one relationship between single-valley-extended solutions with platforms  of Eq.(\ref{FKS2}) and even solutions of Eq.(\ref{FKS}).

\begin{thm}\label{T3}
Suppose that $f$ is a single-valley-extended solutions with platform  of (\ref{FKS2}).
Let $[u,v]$ be the interval consisting of all minimum points of $f$ and $f(\epsilon^2)<u$.
Then $g(x)={\rm sgn}(v-|x|)f(|x|)$ for $x\in [-1,1]$ is an even solution of (\ref{FKS}).
Conversely, if $g$ is an even solution of (\ref{FKS}), then $f(x)=|g(x)|$ for $x\in [0,1]$ is a solution of (\ref{FKS2}).
\end{thm}
\par
\emph{Proof}. We prove the first part.
It follows from Lemma \ref{P4} that
$$
f([0,\epsilon^2u])=[v,1],
f([\epsilon^2u,\epsilon^2v])=[u,v],
f([\epsilon^2v,\epsilon^2])=[\epsilon^2,u].
$$

If $y:=|x|\in [v,1]$, then one has
\begin{eqnarray*}
g(x)&=&{\rm sgn}(v-|x|)f(|x|)=-f(y),\\
\epsilon^2 y&\in& [\epsilon^2v,\epsilon^2]\subset [0, v],\\
f(\lambda y)&\in& [\epsilon^2,u]\subset [0, v].
\end{eqnarray*}
Thus $f(f(\epsilon y)) =g(g(\epsilon^2 x))$. It follows from $f(y)=\epsilon^{-1} f(f(\epsilon^2 y))$
that $g(x)=-\epsilon^{-1} g(g(\epsilon^2 x))$.

If $y:=|x|\in [0,u]$, then one has
\begin{eqnarray*}
g(x)&=&{\rm sgn}(v-|x|)f(|x|)=f(y),\\
\epsilon^2 y &\in& [0,\epsilon^2u]\subset [0, u],\\
f(\epsilon^2 y)&\in& [v,1].
\end{eqnarray*}
Thus $f(f(\epsilon^2 y)) =-g(g(\epsilon^2 x))$. It follows from $f(y)=\epsilon^{-1} f(f(\epsilon^2 y))$
that $g(x)=-\epsilon^{-1} g(g(\epsilon^2 x))$.

If $y:=|x|\in [u,v]$, then one can check
that $g(x)=-\epsilon^{-1} g(g(\epsilon^2 x))=0$.

On the other hand, it is clear that $g(0)=1$ and $g(-x)=g(x)$. Therefore $g$ is an even solution of (\ref{FKS}).

Next we prove the second part.
If $g$ is an even solution of (\ref{FKS}), then
$$
g(x)=-\epsilon^{-1} g(g(\epsilon^2 x))=-\epsilon^{-1}g(|g(\epsilon^2  x)|).
$$
Thus $|g(x)|=\epsilon^{-1} |g(|g(\epsilon^2 x)|)|$.
Therefore $f(x)=|g(x)|$ for $x\in [0,1]$ is a solution of (\ref{FKS2}).
\hfill$\Box$




\end{document}